\def\sqr#1#2{{\vcenter{\hrule height.#2pt
      \hbox{\vrule width.#2pt height#1pt \kern#1pt
         \vrule width.#2pt}
       \hrule height.#2pt}}}
\def\eop{\mathchoice\sqr34\sqr34\sqr{2.1}3\sqr{1.5}3}

\def\N{\hbox{\rm I\hskip -0.14em N}}

\def\mod{\hbox{\rm mod}\hskip 0.5em}
\def\th #1 #2. #3\par{\medbreak{\bf#1 #2.
\enspace}{\sl#3}\par\medbreak}
\documentclass{llncs1}
\usepackage{setspace}
\usepackage{makeidx}  % allows for indexgeneration
\usepackage{graphicx}
\begin{document}
\frontmatter          % for the preliminaries
\pagestyle{headings}  % switches on printing of running heads

\title{On the rank of elliptic curves}
\titlerunning{On the rank of elliptic curves}  % abbreviated title (for running head)
\author{Jorma Jormakka}
\authorrunning{Jorma Jormakka}   % abbreviated author list (for running head)
%
%%%% modified list of authors for the TOC (add the affiliations)
\tocauthor{Jorma Jormakka}

\institute{ Contact by: \email{jorma.o.jormakka@gmail.com}}

\maketitle              % typeset the title of the contribution

\begin{abstract}
The paper proves that the Birch and 
Swinnerton-Dyer conjecture is false. 
\end{abstract}

\begin {keywordname}
Elliptic curves, Euler product, Birch and Swinnerton-Dyer conjecture.
\end{keywordname}

\section{Introduction}

Let 
$P=\{p_1,p_2,\dots| p_j {\rm \ is \ a \ prime, \ } p_{j+1}>p_j>1,j\ge 1\}$
be the set of all primes larger than one.
In [1] an elliptic curve $C$ over the field of rational numbers $Q$ 
is a curve defined by the Weierstrass equation
$$y^2=x^3+ax+b$$
where $a, b\in Z$ and $x,y\in Q$. The discriminant of the cubic equation
is $\Delta=-16(4a^3-27b^2)\not=0$. 
Let $N_{p_j}$ denote the number of solutions to
$y^2=x^3+ax+b \ \mod p_j$ 
and let $a_{p_j}=p_j-N_{p_j}$.

The incomplete L-function of the curve $C$ is 
$$L(C,s)=\prod_{j\in A_C}(1-a_{p_j}p_j^{-s}+p_j^{1-2s})^{-1},\eqno(1)$$
where
$$A_C=\{j\in \N, j>0, p_j \ does \ not \ divide \ \Delta \}.$$
The Euler product (1) converges absolutely at least if $Re\{s\}>2$ because 
$|a_p|\le 2p$. This upper bound for $<a_p|$ is obvious since $x$ takes $p$ 
values and $y$ can take two values for each $x$.
Hasse's statistical bound $|a_p|\le 2p^{1\over 2}$ improves the area of
absolute convergence to $Re\{s\}>3/2$ and [1] gives this area. The problem statement [1] tells that $L(C,s)$ has a holomorphic continuation to the whole 
complex plane, thus it does not have poles.

The Birch and Swinnerton-Dyer conjecture says that the Taylor expansion
of $L(C,s)$ at $s=1$ has the form
$$L(C,s)=c(s-1)^r+ \ higher \ order \ terms\eqno(2)$$
with $c\not=0$ and $r$ the rank of $C$.
The rank of an elliptic curve is defined as the rank of the group of solutions
in the rational numbers. The number $r$ in the Taylor expansion of $L(C,s)$
is called the algebraic rank of the curve. 
The conjecture is thus that the rank and the algebraic rank are equal. 

Let $p>2$ be prime, $Z_p$ the cyclic group of integers
modulo $p$, and $Z_p^*=\{1,\dots ,p-1\}$.
The set of quadratic residues modulo $p$ is the set
$$QR_p=\{ x\in Z_p^* | \exists y\in Z_p^* \ {\rm such} \ {\rm that}\\
 \ y^2\equiv x \ (\mod p)\}$$
and the set of nonresidues modulo $p$ is 
$$QNR_p=\{ x\in Z_p^* | x\not\in QR_p\}.$$
If $g$ is a primitive root of $Z_p^*$, then
$$Z_p^*=\{ g^0,g^1,\dots ,g^{p-2}\}.$$
The set $QR_p$ is the subset where $g$ has even powers:
$$QR_p=\{ g^0,g^2,\dots \}.$$
Thus, \#$QR_p=$\#$QNR_p$, the sets $QR_p$ and $QNR_p$ have 
equally many elements. If the integer $a$ divides integer $b$ it
is written as $a|b$.  
For brevity, we write $y\equiv x$ as a shorthand of $y\equiv x \ (\mod p)$
when there is no chance of confusion.

There is a recursion formula for deriving rational solutions from
a rational base point $(x,y)$
$$x_{i+1}=S_i^2-2x_i \ , y_{i+1}=y_i+S_i(x_{i+1}-x_i)\eqno(3)$$
$$S_i={a+3x_i^2\over 2y_i}$$
This recursion gives a new rational solution in the following way:
$$y_{i+1}^2=y_i^2+2y_iS_i(x_{i+1}-x_i)+S_i^2(x_{i+1}-x_i)^2$$
$$=x_i^3+ax_i+b+2y_iS_i(x_{i+1}-x_i)+S_i^2(x_{i+1}-x_i)^2$$
$$=x_{i+1}^3+ax_{i+1}+b$$
yielding
$$x_{i+1}^2+x_{i+1}x_i+x_i^2+a=2y_iS_i+(x_{i+1}+2x_i)(x_{i+1}-x_i)$$
which gives
$$3x_i^2+a=2y_iS_i.$$
The recursion may end or it may generate an infinite number of 
rational solutions. An example where the recursion ends is the following:

We define the curve $C_1$ by a Weiestrass form with $a=33$ and $b=-26$. As
the base point we take $x_0=3$ and $y_0=10$. The recursion (3) shows
that $x_1=x_0$, $y_1=y_0$. It follows that recursion generates
only one solution $(3,10)$. The curve $C_1$ is a special case of
$$x_0=3s^2\hskip 2em y_0=9s^2\pm s\eqno(4)$$
$$a=27s^4\pm 6s^2\hskip 2em b=s^2-27s^6$$
with $s=1$. For every nonzero integer value $s$ the solution (4) gives 
$x_1=x_0$, $y_1=y_0$. These solutions 
are found by setting $x_1=x_0$ in (3).

The rank of an elliptic curve is the number of independent base points
from which the recursion derives an infinite number of rational solutions.
For $C_1$ the recursion gives only finitely many points, but 
for that special elliptic curve there may be other 
base points that give infinitely many different points. An example of
an elliptic curve having infinitely many rational solutions is
$y^2=x^3-5^2x$. This is known since $5$ is a congruent number.
If $d$ is a noncongruent number, such as $r^2$ for any integer $r$ there
are only three solutions: $(0,0)$ and $(\pm d,0)$.

The recursion formula (3) has a corresponding operation in 
integers modulo $p$ in the form
$$x_{i+1}\equiv S_i^2-2x_i \ (\mod p)\eqno(5)$$
$$S_i\equiv (s+3x_i^2)(2t_i)^{-1} \ (\mod p)$$
$$t_i\equiv x^3+ax+b \ (\mod p).$$

If $(x_i,y_1)$ is a solution in $Z_p^*$ then the recursion formula in 
$Z_p$ gives another solution $(x_{i+1},y_{i+1})$, $x_{i+1},y_{i+1}\in Z_p^*$,
 where 
$$y_{i+1}^2\equiv t_{i+1} \ (\mod p).$$
The operation also takes a pair $(x_i,t_i)$ where $t_i\in QNR_p$ into
a pair $(x_{i+1},t_{i+1})$ where $t_{i+1}\in QNR_p$. Iterating the operation
gives classes of pairs $(x_i,y_i)$. If there
is a solution in $Q$, then all of the iterated solutions map
to the same set of $(x_i,y_i)$ in $Z_p$. 

The claim that the Birch and Swinnerton-Dyer conjecture should hold seems
to be based on the idea that an infinite number of solutions in rationals
for an elliptic curve $C$ would give more solutions in the modular case.
This is a very strange idea because there are very many solutions for a 
modular equation e.g. in knapsack problems and it is very difficult to find
integer solutions to knapsack problems. The modular problem and the integer
problem are quite different issues. The same should be the case with the 
modular problem and the rational problem in elliptic curves. The natural
expectation is that these problems are very different and one does not
give information of the other.    

Two elliptic curves over rationals are 
known to have very high rank (one exactly or rank 20 and the other of 
rank at least 28). They are of the form
$$y^2+xy+y=x^3-x^2-b$$
where $b\in \N$. Writing this curve in the Weierstrass form gives
$$y_1^3=x_1^3+(-{49\over 48})x_1+(-{2149\over 576}-b)$$
$$y_1=y+{1\over 2}(x+1)\hskip 2em x_1=x-{5\over 12}$$
and in the form where coefficients are integers is
$$y_2^3=9x_1^3-147x_1-2149-576b$$
$$y_2=24y+12x+12\hskip 2em x_2=4x-{5\over 3}$$
As $a$ and $b$ in the Weierstrass form are not integers in these elliptic 
curves, they are not elliptic curves considered in [1] and in this paper.

\section{Calculation of $a_{p_j}$ for $y^2=x^3-d^2x$}

\th Lemma 1. 
Let $p>2$ be prime and $a$ an integer. Assume $-1\in QNR_p$ and 
$a\not\equiv 0 \ (\mod p)$. 
The number $N_p$ of solutions to the modular equation 
$$y^2\equiv x^3+ax \ (\mod p) \eqno(6)$$
is $N_p=p.$

\proof
Let
$$A=\{ x\in Z^*_p | t_1\in QR_p, t_1\equiv x(x^2+a) \ (\mod p) \},$$
$$B=\{ x\in Z^*_p | t_1\in QNR_p, t_1\equiv x(x^2+a) \ (\mod p) \},$$
and $m_1=$\#$A$, $m_2=$\#$B$. 
We can write
$$A_1=\{ x=1,\dots, {p-1\over 2}| t_1\in QR_p, t_1\equiv x(x^2+a) \ (\mod p) \},$$
$$A_2=\{ x={p+1\over 2},\dots ,p-1| t_1\in QR_p, t_1\equiv x(x^2+a) \ (\mod p) \},$$
$$B_1=\{ x=1,\dots, {p-1\over 2}| t_1\in QNR_p, t_1\equiv x(x^2+a) \ (\mod p) \},$$
$$B_2=\{ x={p+1\over 2},\dots ,p-1| t_1\in QNR_p, t_1\equiv x(x^2+a) \ (\mod p) \},$$
and $m_{1,i}=$\#$A_i$, $m_{2,i}=$\#$B_i$, $i=1,2$. 
The sets $A_1$ and $A_2$ are disjoint and $A=A_1\cup A_2$.
Similarly, the sets $B_1$ and $B_2$ are disjoint and $B=B_1\cup B_2$.
Calculating
$$A_2=\{ -x=-{p+1\over 2},\dots ,-p-1| t_1\in QR_p, t_1\equiv x(x^2+a) \ (\mod p) \}$$
$$=\{ -x=p-{p+1\over 2},\dots ,p-p+1| t'_1\in QR_p $$
$$\hskip 2em t'_1=p-t_1\equiv (-x)((-x)^2+a) \ (\mod p) \}$$
$$=\{ -x=1,\dots, {p-1\over 2}| t'_1\in QR_p, t'_1\equiv (-x)((-x)^2+a) \ (\mod p) \}.$$
If it were true that $-1\in QR_p$, then there would exists 
$\epsilon$ such that $-1\equiv \epsilon^2$. Then for any $y$ holds 
$-y^2\equiv (\epsilon y)^2\in QR_p$. But as we require that
$-1\in QNR_p$ it is not possible that $-y^2\equiv h^2$ for any $h$ because 
if it is $-1\equiv (y^{-1}h)^2\in QR_p$. Thus, $-y^2\in QNR_p$ for every $y$. 
Therefore
$$A_2=\{ x'=1,\dots, {p-1\over 2}| t'_1\in QNR_p, t'_1\equiv x'(x'^2+a) \ (\mod p) \}=B_1.$$
Similarly, $A_1=B_2$. It follows that
$$m_1=m_{1,1}+m_{1,2}=m_{1,1}+m_{2,1},$$
$$m_2=m_{2,1}+m_{2,2}=m_{2,1}+m_{1,1}.$$
Thus, $m_1=m_2$. Let $a\in QR_p$. Then there are two values $x\in Z^*_p$ 
that yield $t_1\equiv 0 \ (\mod p)$. Therefore
$$m_1+m_2=p-3 \ \Rightarrow m_1={p-3\over 2}.$$
Every $x\in A$ yields two solutions $y, p-y$ to (6). Every $x$ giving
$t_1\equiv 0 \ (\mod p)$ yields one solution $y=0$ to (6). 
The number of solutions is
$$N_p=2{p-3\over 2}+3=p.$$
If $a\in QNR_p$ then $m_1+m_2=p-1$ and 
$$N_p=2{p-1\over 2}+1=p.$$
The lemma is proved. $\eop$

Let us give an example of Lemma 1. Let $d=1$ and $p=11$. Then
$QR_{11}=\{1,3,4,5,9\}$. When $x$ ranges from $0$ to $10$ the values 
of $x(x^2-1)$ give the sequence $0,0,6,2,5,10,1,6,9,5,0$. Removing zeros
from this sequnce as they are neither in $QR_p$ nor in $QNR_p$ we notice
that $-6\equiv 5$. Because $-1\in QNR_p$ we have $6\in QNR_p$ and $-5\in QR_p$.
Likewise $-2\equiv 9$, so $2\in QNR_p$ and $9\in QR_p$; 
$-10\equiv 1$. The same is with $-1\equiv 10\in QNR_p$ 
and $1\in QR_p$. We get $2(p-3)/2=8$ 
solutions: $(4,4)$, $(4,7)$, $(6,1)$, $(6,10)$, $(7,3)$, $(7,8)$,
$(9,4)$, $(9,7)$, that is, for each $x$ there are two $y$ values. 
Additionally we have the zeros. They give three solutions $(1,0)$
and $(10,0)$ from $x^2-1\equiv 0$ and $(0,0)$ is a solution. Together
there are $11=p$ solutions.

\th Lemma 2. 
Let $p>2$ be prime. The number of solutions $y^2$
to the equation
$$y^2-c\equiv x^2 \ (\mod p)\eqno(7)$$
satisfying $y^2,x^2\in Z_p^*$ is 

${p-5\over 4}$ if $-1\in QR_p$ and $c\in QR_p$,

${p-3\over 4}$ if $-1\in QNR_p$,

${p-1\over 4}$ if $-1\in QR_p$ and $c\in QNR_p$.

\proof

Let us assume that (7) holds. Thus there exists $z\in Z_p^*$ such that 
the modular equation
$$y^2-x^2=(y-x)(y+x)\equiv c$$
can be written as
$$y-x\equiv z  \ , \ y+x\equiv z^{-1}c.$$
Then 
$$y\equiv 2^{-1}x^{-1}(z^2+c) \ , \ x\equiv 2^{-1}z^{-1}(z^2-c).$$
Let $\pm \epsilon$ denote the two roots of $z^2\equiv -1$ if $-1\in QR_p$.
If $-1\in QNR_p$ there are no such roots.  

If $c\in QNR_p$ and $-1\in QR_p$ there are no solutions to the equations
$$z^2\equiv c\hskip 2em,\hskip 2em (\epsilon z)^2\equiv -c.\eqno(8)$$ 
In this case we let $z$ range over the $p-1$ numbers in $Z_p^*$ in the
equation for $y$.
If two values $z_1$ and $z_2$ give the same $y$, then
$$z_1^{-1}(z_1^2+c)\equiv z_2^{-1}(z_2^2+c)$$
i.e., 
$$z_1+cz_1^{-1}\equiv z_2+cz_2^{-1},$$
$$z_1-z_2\equiv c(z_2^{-1}-z_1^{-1}).$$
Multiplying by $z_1z_2$
$$z_1z_2(z_1-z_2)\equiv z_1z_2c(z_2^{-1}-z_1^{-1})\equiv c(z_1-z_2)$$
and $z_1z_2\equiv c$, i.e, $z_2\equiv cz_1^{-1}$. 
When $z$ ranges over
all values in $Z_p^*$ the number $y$ gets all values it can get 
and exactly two values $z$ map to the same $y$. 
The number of different $y$ is therefore ${p-1\over 2}$.

If some value of $z$ gives $y$, another value of $z$ gives $-y$. 
As $\pm y$ yield the same $y^2$ the number of different $y^2$ 
is half of the numbers of $y$, that is, ${p-1\over 4}$. 

If $c\in QR_p$ and $-1\in QR_p$ then there are two solutions $z$ to
both of the equations in (8). These four values of $z$ are all different.
Removing them gives $p-5$ values for the range of $z$. 
The number of different values $y^2$ is ${p-5\over 4}$.  

If $c\in QR_p$ and $-1\in QNR_p$ there are two solutions for
$z^2\equiv c$ but no solutions to $z^2\equiv -c$. 
The number of different $y^2$ is ${p-3\over 4}$.  

If $c\in QNR_p$ and $-1\in QNR_p$ there are no solutions for
$z^2\equiv c$ but two solutions to $z^2\equiv -c$. 
The number of different $y^2$ is ${p-3\over 4}$.$\eop$

\th Lemma 3. 
Let $p>2$ be prime and $a$ an integer. Let $-1\in QR_p$,
$a\not\equiv 0 \ (\mod p)$ and $g$ a primitive root of $Z_p^*$. 
The number $N_p$ of solutions to the modular equation 
$$y^2\equiv x^3+ax \ (\mod p) \eqno(9)$$
is 

$N_p=8n_1+7$ if $-a\equiv g^{2i}$ and $i$ is even,

$N_p=2p-8n_1-7$ if $-a\equiv g^{2i}$ and $i$ is odd,

$N_p=8n_g+3$ if $-a\equiv g^{2i+1}$ and $i$ is even,

$N_p=2p-8n_g-3$ if $-a\equiv g^{2i+1}$ and $i$ is odd.

Here $n_c$ is the number of solutions $y^4\in Z_p^*$ 
yielding $y^4-c\in QR_p$, $c=1$ or $c=g$.

\proof
Let 
$$B=\{ x'\in Z^*_p | t'\in QR_p, t'\equiv x'^3+ax' \ (\mod p) \}.\eqno(10)$$

If $-a\equiv g^{2i}$ we insert $t\equiv g^{-3i}t'$ and $x\equiv g^{-1}x'$. 
Then $t'\equiv x'^3+ax'$ changes to $g^{3i}\equiv g^{3i}x^3-g^{2i}g^ix$,
i.e., to $t\equiv x^3-x$. We reduced $-a$ to $c=1$. 

If $-a\equiv g^{2i+1}$ we insert the same $t$ and $x$ as above. Then
$t'\equiv x'^3+ax'$ changes to $g^{3i}\equiv g^{3i}x^3-g^{2i+1}g^ix$,
i.e., to $t\equiv x^3-gx$. We reduced $-a$ to $c=g$.

We write both of these cases as $t\equiv x^3-cx$ where $c=1$ if 
$-a=g^{2i}$ and $c=g$ if $-a=g^{2i+1}$.   
 
Let
$$A= \{ x\in Z^*_p | t_1\in QR_p, t_1\equiv x^3-cx \ (\mod p) \}\eqno(11)$$
$$A'=\{ x\in Z^*_p | t_1\in QNR_p, t_1\equiv x^3-cx \ (\mod p) \}.$$

If $i$ is even then $g^i$ is in $QR_p$ and in the substitution $t=g^{3i}t'$
holds: if $t\in QR_p$ then $t'\in QR_p$. 
If $i$ is odd, then $t\in QR_p$ implies 
that $t'\in QNR_p$. Thus, for even $i$ $B=A$ while for odd $i$ $B=A'$.

Let us write the sets $A$ and $A'$ differently
$$A=\{g^k|g^k(g^{2k}-c) \in QR_p , k=0,\dots, p-2 \} \eqno(12)$$
$$A'=\{g^k|g^k(g^{2k}-c) \in QNR_p , k=0,\dots, p-2 \}$$
and let us divide them into subsets of even and odd indices of $k$
$$A_1=\{g^{2k}|g^{2k}(g^{4k}-c) \in QR_p, k=0,\dots, {p-3\over 2} \}$$
$$A_2=\{g^{2k}|g^{2k+1}(g^{2(2k+c)}-c) \in QR_p, k=0,\dots, {p-3\over 2} \}$$
$$A'_1=\{g^{2k+1}|g^{2k}(g^{4k}-c) \in QNR_p, k=0,\dots, {p-3\over 2} \}$$
$$A'_2=\{g^{2k+1}|g^{2k+1}(g^{2(2k+c)}-c) \in QNR_p, k=0,\dots, {p-3\over 2} \}.$$
Then $A=A_1\cup A_2$, \#$A=$\#$A_1+$\#$A_2$ and
$A'=A'_1\cup A'_2$, \#$A'=$\#$A'_1+$\#$A'_2$.
 
We also define sets that do not have the $x=g^k$ term in $t=x(x^2-c)$. 
$$C=\{g^{2k}|g^{2k}-c \in QR_p, k=0,\dots, {p-3\over 2} \}\eqno(13)$$
$$C'=\{g^{2k}|g^{2k}-c \in QNR_p, k=0,\dots, {p-3\over 2} \}$$
and divide these sets into subsets where a set with a running index 
$2k$ is divided into two sets with running indices $4k$ and $4k+2$:
$$C_1=\{g^{4k}| g^{4k}-c \in QR_p, g^{4k}\le {p-1\over 2}, 
k=0,\dots, {p-3\over 2} \}\eqno(14)$$
$$C_2=\{g^{4k+2}| g^{4k+2}-c \in QR_p, g^{4k+2}\le {p-1\over 2}, 
k=0,\dots, {p-3\over 2} \},$$
$$C'_1=\{g^{4k}| g^{4k}-c \in QNR_p , g^{4k}\le {p-1\over 2}, 
k=0,\dots, {p-3\over 2} \},$$
$$C'_2=\{g^{4k+2}| g^{4k+2}-c \in QNR_p, g^{4k+2}\le {p-1\over 2}, 
k=0,\dots, {p-3\over 2} \}.$$
The rule $g^{4k}\le {p-1\over 2}$ and $g^{4k+2}\le {p-1\over 2}$ 
removes half of the values of the running index.
Then $C=C_1\cup C_2$, \#$C=$\#$C_1+$\#$C_2$ and
$C'=C'_1\cup C'_2$, \#$C'=$\#$C'_1+$\#$C'_2$.

The idea is to map the solutions of $t\equiv g^k(g^{2k}-c)$ bijectively to 
solutions of $t'\equiv g^{2k}-c$.
Clearly, if $g^k\in QNR_p$ multiplying with it changes $t'\in QR_p$ to
$t\in QNR_p$ and if $g^k\in QR_p$ multiplying by it does not change the set.
This is why we divided the sets to $A_i$, $A'_i$, $i=1,2$. 
In $i=2$ sets $g^k\in QNR_p$, 
so if an element of $C'_2$ is multiplied by $g^k$ we get an element of
$A_2$. Likewise, $C_2$ and $A'_2$ correspond to each other.    
  
The following relations hold

\centerline {\#$A=$\#$2C_1+$\#$2C'_2$}

\centerline {\#$A'=$\#$2C'_1+$\#$2C_2$}

\centerline {\#$C_2=$\#$C-$\#$C_1$}

\centerline {\#$C'_2=$\#$C'-$\#$C'_1$.}

Solving \#$A$ yields

\centerline {\#$A=2$\#$C_1+2$\#$C'_2$}

\centerline {$=2$\#$C_1+2$\#$C'-2$\#$C'_1$.}

The value $a$ is used in the proof of this lemma in two places only. 
One is in (10): if $-a\equiv g^{2i}$ or $-a\equiv g^{2i+1}$ 
and the index $i$ is even, then $B=A$. If $i$ is odd, then $B=A'$.
The other place is in Lemma 2 where the numbers of solutions in the
different cases depend on if whether $-a\equiv g^{2i}$ or 
$-a\equiv g^{2i+1}$, i.e., if $c=1$ or $c=g$.

Case 1: $-a=g^{2i}$, $i$ even. Then $c=1$ and the relation 

\centerline {\#$C'={p-3\over 2}-$\#$C$.}

holds. In this relation we have counted the values of $k$ in
$C\cup C'$ and excluded the one value of $k$ that gives 
$g^{2k}-1\equiv 0 \ (\mod p)$ because $0\not \in QR_p\cup QNR_p$.
Thus, the number of valid indices $k$ is one less than the number
${p-1\over 2}$ of indices $k=0,\dots,{p-3\over 2}$ in (13). 
The correct value of valid indices ${p-3\over 2}$.

Counting indices $k$ in  $C_1\cup C'_1$ in (14) gives 

\centerline {\#$C'_1={p-5\over 4}-$\#$C_1.$}

In this relation we have excluded the solution to
$g^{2k}-c\equiv 0 \ (\mod p)$. In $C_1$ and $C'_1$ the counted element is 
not the number of indices $k$. It is the number of values $g^{4k}$.  
This number must be reduced by one. The result follows as
${p-1\over 4}-1={p-5\over 4}$.

Calculating \#$A$ gives

\centerline {\#$A=2$\#$C_1+p-3-2$\#$C-{p-5\over 2}+2$\#$C_1$}

\centerline {$=4$\#$C_1-2$\#$C+{p-1\over 2}.$}

Writing \#$C_1=n_1$ and inserting from Lemma 2 the case
$-1\in QR_p$ and $c\in QR_p$ where \#$C={p-5\over 4}$ yields

\centerline {\#A$=4n_1-{p-5\over 2}+{p-1\over 2}=4n_1+2.$}

In Case 1 holds $N_p=2$\#$A+3$ because if there is a solution
$y^2\equiv x(x^2-1)$, then it is satisfied by two $y$ values, $\pm y$, 
and there are three solutions where $y\equiv 0$, namely $x\equiv 0$,
$x^2\equiv \pm 1$. Thus $N_p=8n_1+7$.

Case 2: $-a=g^{2i}$, $i$ odd. Then $B=A'$. Thus

\centerline {\#$A'=2$\#$C'_1+2$\#$C_2$}

\centerline {$={p-5\over 2}-2$\#$C_1+2$\#$C-2$\#$C_1$}

\centerline {$={p-5\over 2}-4n_1+{p-5\over 2}=p-4n_1-5$}

In Case 2 $N_p=2$\#$A'+3$, thus $N_p=2p-8n_1-7$.

Case 3: $-a=g^{2i+1}$, $i$ even. The differences to Case 1 are

\centerline {$N_p=2$\#$A+1$} 

\centerline {\#$C'={p-1\over 2}-$\#$C$}

\centerline {\#$C'_1={p-3\over 4}-$\#$C_1$}

\centerline {\#$C={p-3\over 4}$}

because $c=g$ and $g^{2k}-g\equiv 0$ is not possible.  

We denote \#$C_1=n_g$ and insert from Lemma 2 the case
$-1\in QR_p$ and $c\in QNR_p$ where \#$C={p-1\over 4}$.
Making these changes to the calculation of Case 1 gives
$N_p=8n_g+3$. 

Case 4: $-a=g^{2i+1}$, $i$ odd. 
Analogically with Cases 2 and 3 we get
$N_p=2p-8n_g-3$.

The definition of $C_1$ is a bit complicated as the running index $k$ 
loops over twice as many indices than are needed and the set has 
a test to discard half of the values values $k$ because in this
way $C$ is clearly the union of $C_1$ and $C_2$.  
It is good to notice that the set $C_1$ has as many members as the set
$$\{y^4|y^4-c\in QR_p\}$$
where $c=1$ for Cases 1 and 2 and $c=g$ for Cases 3 and 4.  

The proof of the lemma is completed. $\eop$

Let us look at an example of Lemma 2. 
Let $a=-d^2$ for $d=1$ in (9) and $p=13$. As $g$ we choose $2$, which 
is a primitive root for $Z_{13}^*$.
Then $1\equiv 2^0$, $2\equiv 2^1$, 
$3\equiv 2^4$, $4\equiv 2^2$, $5\equiv 2^9$, $6\equiv 2^5$,
$7\equiv 2^{11}$, $8\equiv 2^3$, $9\equiv 2^8$, $10\equiv 2^{10}$,
$11\equiv 2^7$ and $12\equiv 2^6$. The sets are
$$A=\{2^3,2^9\}\hskip 2em A'=\{2^1,2^2,2^4,2^5,2^7,2^8,2^{10},2^{11}\}$$
$$A_1=\emptyset\hskip 2em A_2=\{2^3,2^9\}$$
$$A'_1=\{2^2,2^4,2^8,2^{10}\}\hskip 2em A'_2=\{2^1,2^5,2^7,2^{11}\}$$
$$C=\{2^1,2^5\}\hskip 2em C'=\{2^2,2^3,2^4,2^9\}$$
$$C_1=\emptyset\hskip 2em C_2=\{2^1,2^5\}$$
$$C'_1=\{2^2,2^4\}\hskip 2em C'_2=\{2^9\}$$                
There is a direct correspondence between $C_1$ and the first half of $A_1$,
as there is between $C'_1$ and the second half of $A'_1$. This is because if 
$g^{2k}(g^{4k}-1) \in QR_p$ 
then $g^{4k}-1 \in QR_p$ 
and if $g^{4k}-1 \in QR_p$ then $\pm g^{2k}(g^{4k}-1) \in QR_p$ since
$-1\in QR_p$.  
There is also a direct  
correspondence between $C'_2$ and the first half of $A_2$,
as there is between $C_2$ and $A'_2$. This is because if 
$g^{2k+1}(g^{2(2k+1)}-1) \in QR_p$ 
then $g^{2(2k+1)}-1 \in QNR_p$ as $g\in QNR_p$, 
and if $g^{2(2k+1)}-1 \in QNR_p$ then 
$\pm g^{2k+1}(g^{2(2k+1)}-1) \in QR_p$ since $-1\in QR_p$. This gives
the relations between the sizes of the sets. 

The case of Lemma 1 covers half of all $p$ because of Lemma 4.

\th Lemma 4.
The following statements hold:

(i) $-1\in QR_p$ if and only if $4|(p-1)$

(ii) The number of $p<N$ such that $4|(p-1)$ approaches half when $N$ 
grows to infinity.

\proof 
Let $g$ be a primitive root of $Z_p^*$. If $4|(p-1)$, then
$a\equiv g^{p-1\over 4}$ is in $Z_p^*$ and $a^2\equiv -1$.
If $-1\in QR_p$, then $-1\equiv g^{2i}$ for some $i$, $0\le i\le p-2$.  
Since $-1\not\equiv 1$ holds 
$2i\not\equiv 0 \ (\mod (p-1))$. Thus $2i\not=0$ and $2i\not= p-1$.
As $(-1)^2\equiv 1\equiv g^{p-1}\equiv g^{4i}$ holds
$4i=k(p-1)$ for some $k$ where $k$ has the possible values $1,2,3$.
If $k=2$, then $-1\equiv g^{2i}\equiv g^{p-1}\equiv 1$, which is impossible.
Thus, $k\in \{1,3\}$. Then $gcd(4,k)=1$ and therefore $4|(p-1)$. This 
proves the claim (i).

Claim (ii) is shown true by considering the Sieve of Eratosthenes. 
In this algorithm primes are found by reserving a memory vector for all 
numbers and marking the place of $1$ as full and all other places empty 
at the beginning.
On each step the first unmarked place is taken as the next prime $p$. The
place of $p$ is marked and all multiples of $p$ are marked. In this algorithm
the first step takes $p=2$ and marks all multiples of $2$.
The unmarked numbers are all odd. The next prime is $p=3$,
the first unmarked number. All multiples of $3$ are marked. 
The numbers that are marked for $p$, i.e., multiples of $p$, are
are all odd and equally distributed modulo $4$. Consequently, the 
numbers that remain unmarked are all odd and equally distributed
modulo $4$. This continues in each step, thus the 
numbers that remain unmarked are all odd and equally distributed
between $1$ (mod $4$) and $3$ (mod $4$).

In each step the first unmarked number is the next prime $p$. It is
selected as the smallest number in a set of unmarked numbers that are
always odd and distributed equally between two sets
$1$ (mod $4$) and $3$ (mod $4$).
The next prime $p$ has half a chance in belonging to either set. 
The number $p-1$ is always even and if $p\equiv 1$
modulo $4$, then $4|(p-1)$. This is so in half of the cases when $N$ 
approaches infinity. $\eop$

The numbers $n_1$ and $n_2$ in Lemma 3 are not easily evaluated 
for statistical purposes as we did not check how many numbers $x$
give the same $t\equiv x(x^2+a)$. This is done in Lemma 5.

\th Lemma 5. Let $p>2$ be a prime and $-1\in QR_p$.
If $c=1$ the solutions are divided into singlets and multiplets. In
the set of singlets each $x$ is mapped to a unique $y^2$. In multiplets
three, in maximum size cases two, values of $x$ are mapped to the same 
$y^2$. 

(i) There is a running index $t$ such that in singlets $t$ runs
from $1$ to $p-1$ and two values of $t$ map to the same $x$. 

(ii) In multiplets there is a running index $t$ such that $t$ runs 
from $1$ to $p-1$ and two values of $t$ map to the same $x$. 
For $p-m$ values of $t$ three values of $x$ map to the same $y^2$ and
three values of $-x$ map to the same $-y^2\in QR_p$. Thus, twelve indices $t$ 
map together either to $QR_p$ or to $QNR_p$. 
That is, two $t$ map to the same $x$. Six $t$ map to the same $(h,x)$ where
$h\equiv x(x^2-c)$. If $h\in QR_p$, then $-h\in QR_p$ and twelve values of
$t$ giving $(h,x)$ or $(-h,-x)$ are all in $QR_p$. If $h\in QRP_p$, then
twelve values of $t$ map to $(h,x)$ or $(-h,-x)$ are all in $QNR_p$.
The value of $m$ is as follows:

Case $c=1$:
$$m=13 \ {\rm if} \ 3\in QR_p \ {\rm and} \ 3^{-1}+1\in QR_p$$
$$m=9 \ {\rm if} \ 3\in QNR_p \ {\rm and} \ 3^{-1}+1\in QR_p$$
$$m=9 \ {\rm if} \ 3\in QR_p \ {\rm and} \ 3^{-1}+1\in QNR_p$$
$$m=5 \ {\rm if} \ 3\in QNR_p \ {\rm and} \ 3^{-1}+1\in QNR_p$$

Case $c=g$:

\centerline {If $3\in QR_p$ and $g+1\in QR_p$ then $m=5$.} 

\centerline {If $3\in QR_p$ and $g+1\in QNR_p$ then $m=1$.}

\centerline {If $3\in QR_p$ then $m=1$.}

(iii) In multiplets there are $m_1$ values of $t$ that map to $x$ such that
two values of $x$ map to the same $h\equiv x(x^2-c)$ 
because one of the three values of $x$ is equal to one of the other two. 
In this case eight values of $t$ map
to $(h,x)$ or $(-h,-x)$, all are either in $QR_p$ or all in $QNR_p$.       
The value of $m_1$ is as follows:

Case $c=1$:
$$m_1=6 \ {\rm if} \ 3\in QR_p \ {\rm and} \ 3^{-1}+1\in QR_p$$
$$m_1=4 \ {\rm if} \ 3\in QNR_p \ {\rm and} \ 3^{-1}+1\in QR_p$$
$$m_1=2 \ {\rm if} \ 3\in QR_p \ {\rm and} \ 3^{-1}+1\in QNR_p$$
$$m_1=0 \ {\rm if} \ 3\in QNR_p \ {\rm and} \ 3^{-1}+1\in QNR_p$$

Case $c=g$:

\centerline {If $3\in QR_p$ and $g+1\in QR_p$ then $m_1=4$.}

\centerline {If $3\in QR_p$ and $g+1\in QNR_p$ then $m_1=0$.} 

\centerline {If $3\in QR_p$ then $m_1=0$.}

\proof
As in Lemma 3 we can reduce the equation $h=x(x^2+a)$ to 
$h=x(x^2-c)$ where $c=1$ or $c=g$, where $g$ is a primitive root
of $Z_p^*$. We will consider only the cases $c=1$ and $c=g$.

Let us first consider how the equation
$$x^2-r^2\equiv y^2$$ 
where $x$ and $y$ belong to $QR_p\cup QNR_p$ and $r\not\equiv 0$ is fixed
can be solved with a running index $t$. We have two ways of solving it. 

In the first way we write $z^2\equiv y^2$ and
$$x^2-r^2\equiv (x+r)(x-r)\equiv z^2.$$
There exist $t$ such that
$$x+r\equiv zt \ , \ x-r\equiv zt^{-1}.$$
Then 
$$x\equiv 2^{-1}z(t+t^{-1}) \ , \ z\equiv 2r(t-t^{-1})^{-1}.$$
Thus
$$x\equiv r(t-t^{-1})^{-1}(t+t^{-1}) \ , \ y\equiv \pm 2r(t-t^{-1})^{-1}.$$
If $x\equiv r(t-t^{-1})^{-1}(t+t^{-1})$ for some $t$, then the equation is 
automatically satisfied. We can count the indices $t$.
If two values $t_1$ and $t_2$ give the same $x$, then
$$(t_1-t_1^{-1})^{-1}(t_1+t_1^{-1})\equiv (t_2-t_2^{-1})^{-1}(t_2+t_2^{-1}).$$
Solving gives $t_2^2\equiv t^2_1$. There are three different values 
of $t$, namely $t_1, t_1^{-1}, -t_1^{-1}$, that could give $x$, but 
checking shows that $t_1^{-1}$ gives $-x$ and only two different values 
of $t$, namely $t_1$ and $-t_1^{-1}$ give $x$. Additionally $x$ and $-x$ 
yield the same $y^2$ and $\pm y$ give the same $y^2$. Thus, eight indices
$t$ yield the same solution. 

The second way to solve the equation is to write it as
$$x^2-z^2\equiv r^2 \ , \ z^2\equiv y^2$$
and find a running index $t$ such that
$$x+z\equiv rt \ , \ x-z\equiv rt^{-1}.$$  
Then $x=2^{-1}r(t+t^{-1})$ and $y=\pm 2^{-1}r(t-t^{-1})$.
If $x\equiv 2^{-1}r(t+t^{-1})$ for some $t$, then the equation is 
automatically satisfied. Two values of $t$ map to the same $x$ because if
$$t_1+t_1^{-1}\equiv t_2+t_2^{-1}$$
then $t_2\equiv t_1^{-1}$. The two values $\pm x$ give the same $x^2$
and two values $\pm y$ give the same $y^2$. Also in this way of solving
eight indices of $t$ give the same solution. The two solution methods give
the same result. 

In both ways we have to exclude indices $t$ for which 
$t-t^{-1}\equiv 0$ or $t+t^{-1}\equiv 0$ because for suc $t$ the number 
$x$ or $y$ is not in $QR_p\cup QNR_p$.

Next, consider the case of
$$x^2+a\equiv b$$
where $a\in QNR_p$ and $b\in QNR_p$. Multiplying by a primitive root $g$
gives the equation 
$$c^2-d^2\equiv (c+d)(c-d)\equiv gx^2$$
where $c^2\equiv gb\in QR_p$, $d^2\equiv ga\in QR_p$ and we can use 
the running index method as
$$c+d\equiv gxt \ , \ c-d\equiv xt^{-1}.$$
If $b\in QNR_p$ and $-a\in QR_p$ we can write $-a=c^2$ and proceed as
$$x^2-c^2\equiv (x+c)(x-c)\equiv b.$$
The running index solution method fails only if $b\in QNR_p$ and $a\in QR_p$ 
but $-1\in QNR_p$.

The second general issue we have to look is how to estimate the
number of solutions to $h\in QR_p$ for $h\equiv x(x^2-c)$.
We can solve the equation
$$z^2\equiv x^2-c$$
by a running index $s$ as 
$$z+x\equiv -cs \ , z-x\equiv s^{-1}$$
i.e., if $x\equiv -2^{-1}(cs+s^{-1})$ for some $s$, then $x^2-c\in QR_p$.
If $x$ is not of that form for any $s$, then $x^2-c\in QNR_p$. 
If ($x\in QR_p$ and $x^2-c\in QR_p$) or ($x\in QNR_p$ and $x^2-c\in QNR_p$), 
then $h\equiv x(x^2-c)\in QR_p$, else $h\in QNR_p$.
The numbers $x$ are produced by a running index $t$ that gives
every number $x$ twice. Thus, $x$ has half of the values in $Z_p^*$,
or possibly a few less if some $t$ must be discarded. The possible values 
that $x$ should have in order that $x^2-c$ is a square $z^2$ are given 
by another running index $s$. Also here some values of $s$ may need to be 
discarded. The running indes $s$ also gives each $x$ twice and can reach
(about) half of the values in $Z_p^*$.   

If these two sets of numbers $x$ can be considered independent, then 
for one fourth of the indices $t$ holds $x\in QR_p$ and $x^2-c\in QR_p$.
For one fourth of the indices $t$  holds $x\in QNR_p$ and $x^2-c\in QNR_p$.
Thus, for half of the indices $t$ holds $h\in QR_p$. If so, then
we can think of indices $t$ as probabilistic trials where half of the time 
hit a success, $h\in QR_p$, and half of the time fail, $h\in QNR_p$.   
We will later formulate this condition as the Statistical Assumption. 

The issue in this statistical method to know is how many of these 
trials can be considered independent.
It means, how many indices $t$ act as a group where all values of $t$ 
in the group give either $h\in QR_p$ or $h\in QNR_p$. 
We alreasy know that at least four
indices $t$ act as a group: each $x$ is given by two values of $t$ and
$\pm x$ map to $\pm h$, which are either both in $QR_p$ or both in $QNR_p$
when $-1\in QR_p$. There can be more indices $t$ acting as a group.

Let us start from the case where more than one $x$ map to the 
same $h\equiv x(x^2-c)$. If so, then
we can find $x_1$ and $x_2$, $x_1-x_2\not\equiv 0$, such that
$$x_1^3-cx_1\equiv x_2^3-cx_2.$$
Then 
$$x_1^2+x_1x^2+x_2^2\equiv c$$
i.e.,
$$(2x_1+x_2)^2+3x_2^2\equiv 4c.\eqno(15)$$
If there are solutions $x_1\not\equiv x_2$ to (15), we say that $x_2$
belongs to the case of multiples. 
If $x_2$ is the only value of $x$ that gives $t$, then we
say that $x_2$ belongs to the case of singletons. We have to look separately
at the cases $c=1$ and $c=g$. 

Case $c=1$. Let us write
$$z\equiv \pm (2x_1+x_2)$$
and 
$$2^2-z^2\equiv 3x_2^2.$$
Then there exists $t$ such that
$$2+z\equiv 3tx_2 \ , \ 2-z\equiv t^{-1}x_2.$$
Solving these equations yields
$$x_2\equiv 4(3t+t^{-1})^{-1} \ , z\equiv 2^{-1}x_2(3t-t^{-1}).$$
The value of $x_{1\pm}$ is then
$$x_{1\pm}=2^{-1}(z-x_2)\equiv 2^{-2}x_2(\pm (3t-t^{-1})-2).$$
We can write
$$h\equiv x_2(x_2^2-c)\equiv 4(3t+t^{-1})^{-2}(3t+t^{-1})^{-1}(16-(3t+t^{-1})^{2}).\eqno(16)$$

Let us do the former general consideration with two running indices $t$ and $s$
explicitly, just for clarity.
We take another running index $s$ that assures that the term 
$16-(3t+t^{-1})^{2}$ 
is a square. As we assume that $-1\in QR_p$ in this lemma, we
can write the square as $-r^2$ for some $r$. Thus
$$16-(3t+t^{-1})^{2}\equiv -r^2\eqno(17)$$
which yields
$$(3t+t^{-1})^2-r^2\equiv 16.$$
There exists $s$ such that
$$3t+t^{-1}+r\equiv 4s \ , \  3t+t^{-1}-r\equiv 4s^{-1}.$$
Solving gives
$$3t+t^{-1}\equiv 2(s+s^{-1}).\eqno (18)$$
If $s$ loops over all numbers in $Z_p^*$ and $t$ satisfies (18), then
the term (17) is a square in (16) and 
$$h\equiv x_2(x_2^2-c)\equiv m^2 (3t+t^{-1})^{-1}$$
where $m^2$ is a square:
$$m^2\equiv 4(3t+t^{-1})^{-2}(16-(3t+t^{-1})^{2})\equiv -4(s+s^{-1})^{-2}(s-s^{-1})^2.$$
Writing $3t+t^{-1}$ as $2(s+s^{-1})$ gives
$$h\equiv x_2(x_2^2-c)\equiv m^2 2(s+s^{-1})^{-1}$$
showing very clearly that $h\in QR_p$ if and only if $2(s+s^{-1})\in QR_p$.

If $16-(3t+t^{-1})^{2}$ is not a square, then $3t+t^{-1}$ belongs to
the numbers that are not reached by $s+s^{-1}$. This set of numbers
is (about) half of all numbers in $Z_p^*$. If $3t+t^{-1}$ is in that set and
if $(3t+t^{-1})\in QNR_p$, then from (9) follows that
$$b\equiv (3t+t^{-1})^{-1}(16-(3t+t^{-1})^{2})$$
is in $QR_p$ and 
$$h\equiv x_2(x_2^2-c)\equiv 4(3t+t^{-1})^{-2}b$$
is in $QR_p$. This explicit calculation agrees with the general consideration
that $h\in QR_p$ for (about) half of the indices $t$, but there is
some probabilistic variation, the trials of $t$ can be seen as probabilistic
trials.  

In the case of multiplets when $c=1$ we may have to discard a few cases 
of the running index $t$: if we do not then some cases may do not give 
three different values of $x$, or any values at all.
As  $-1\in QR_p$ is assumed in the lemma, there exists $\epsilon$ such that
$\epsilon^2\equiv -1$.
If $3\in QR_p$, then $3\equiv \beta^2$ for some $\beta$. 
Then $t$ values $\pm \epsilon\beta$ and $\pm \beta$ need to be discarded.
If $3\in QNR_p$, there are no solutions to $t^2\equiv \pm 3$.  

The number $x_2$ cannot be modulo zero because zero is not in
$QR_p$ or $QNR_p$.
Thus, the roots of $3t+t^{-1}\equiv 0$ must be discarded.
This means that if $3\in QR_p$, two values of $t$ must be discarded. 

There are cases when one of the three $x$ values equals one of the two.
If $x_{1+}\equiv x_{1-}$, then $z\equiv 0$. 
This means that $3t-t^{-1}\equiv 0$, i.e., $t^2\equiv 3^{-1}$.  
For $3\in QNR_p$ this cannot happen,
but for $3\in QR_p$ there are two values of $t$ filling this condition.
If so, we get only two different values of $x$. 

If $x_2\equiv x_{i\pm}$, then $\pm 6\equiv 3t-t^{-1}$, i.e.,
$(t\pm 1)^2\equiv 3^{-1}+1$. If $3^{-1}+1\in QR_p$, there are four
values of $t$ that give only two different values of $x$.

Solutions $x\equiv \pm 1$ give $x^2-c\equiv 0$ and thus $h\equiv 0$.
We have to discard four values of $t$ that give $x\equiv \pm 1$.

Counting the numbers $m$ and $m_1$ gives:
$$p-m=p-1-2-2-4-4=p-13 \ {\rm if} \ 3\in QR_p \ {\rm and} \ 3^{-1}+1\in QR_p$$
$$p-m=p-1-4-4=p-9 \ {\rm if} \ 3\in QNR_p \ {\rm and} \ 3^{-1}+1\in QR_p$$
$$p-m=p-1-4-4=p-9 \ {\rm if} \ 3\in QR_p \ {\rm and} \ 3^{-1}+1\in QNR_p$$
$$p-m=p-1-4=p-5 \ {\rm if} \ 3\in QNR_p \ {\rm and} \ 3^{-1}+1\in QNR_p$$
$$m_1=6 \ {\rm if} \ 3\in QR_p \ {\rm and} \ 3^{-1}+1\in QR_p$$
$$m_1=4 \ {\rm if} \ 3\in QNR_p \ {\rm and} \ 3^{-1}+1\in QR_p$$
$$m_1=2 \ {\rm if} \ 3\in QR_p \ {\rm and} \ 3^{-1}+1\in QNR_p$$
$$m_1=0 \ {\rm if} \ 3\in QNR_p \ {\rm and} \ 3^{-1}+1\in QNR_p$$

For singletons Lemma 3 gives a calculation of $n_1$ where the looping index $t$
goes from $t=1$ to $t=p-1$ and gives each $x$ two times.
Exactly half of the numbers $t$ give multiplets and are first removed. 
The remaining (exact) half of the numbers $t$ give singletons.
The cases to be excluded are already in the multiplet numbers. Thus, 
singletons give exactly $(p-1)/2$ values of $x$. The values $\pm x$ yield 
the same $h$ and two values of $y$ give $h\equiv y^2\in QR_p$. Thus,
eight values of the running index $t$ map as a group. This is shown
in Lemma 3 where $N_p=8p+7$ for Case 1. About half of the groups map to
$h\in QR_p$ and about half to $h\in QNR_p$.

Case $c=g$. For $x_2$ belonging to the case of multiplets (15) gives 
the equation
$$(2x_1+x_2)^2+3x_2^2\equiv 4g.\eqno(19)$$
We split the analysis into two cases:

Case $3\in QNR_p$.
Then we write $z^2\equiv (2x_1+x_2)^2$ and thus 
$$g^{-1}z^2\equiv 2^2-(\eta x_2)^2$$
where $\eta^2\equiv 3g^{-1}$. Then 
$$2+\eta x_2\equiv g^{-1}zt \ , \ 2-\eta x_2\equiv zt^{-1}.$$
Thus $4\equiv z(g^{-1}t+t^{-1})$ and $x_1\equiv 2^{-1}(\pm z-x_2)$, i.e.,
$$x_2\equiv (2\eta)^{-1}4(g^{-1}t+t^{-1})^{-1}(g^{-1}t-t^{-1}).$$
There are two values of $t$ that map to the same $x_2$, namely
$t$ and $-t^{-1}$, solved as before. The third value $t^{-1}$ maps to 
$-x_2$. 

There are no values of $t$ giving $g^{-1}t+t^{-1}\equiv 0$ which would yield 
no $z$ or no $x_2$ and would have to be discarded. 
Thus, all values of $t$ yield $x_2$ and $x_{1\pm}$.  
There are no values of $t$ for which $x_2^2-g\equiv 0$ and 
consequently $h\equiv 0$.
Thus, all values of $t$ yield a valid $h$. 
There are no values of $t$ for which $z\equiv 0$ and 
$x_{1+}\equiv x_{1-}$. All these special cases are missing if $c=g$ and
$3\in QNR_p$, but there is still one special case left:

If $x_2\equiv x_{1\pm}$ then $3x_2\equiv \pm z$ and solving it gives
$$t^2\pm 2\eta 3^{-1}g t-g\equiv 0.$$
Inserting $3g^{-1}\equiv \eta^2$ yields  
$$t^2\pm 2t+1\equiv (t\pm 1)^2 \equiv g+1.$$
If $g+1\in QR_p$ 
there are four indices of $t$ that give solutions
where $x_2\equiv x_{1+}$ or $x_2\equiv x_{1-}$. The other values of
$t$ give three different values of $x$. If $g+1\in QR_p$ we get $m=5$ and 
$m_1=4$. If $g+1\in QNR_p$ then $m=1$, $m_1=0$.

Case $g=1$ and $3\in QR_p$.
Then there exists $\beta$ such that $\beta^2\equiv 3$ and since we assume
that $-1\in QR_p$ there exists $\epsilon$ such that $\epsilon^2\equiv -1$.
Then $-3\equiv (\epsilon\beta)^2$ and (12) can be written as
$$z^2-(\epsilon\beta x_2)^2\equiv 4g$$
where 
$$z\equiv \pm 2x_1+x_2.$$
There exists $t$ such that
$$z+\epsilon\beta x_2\equiv 2gt \ , z-\epsilon\beta x_2\equiv 2t^{-1}.$$
Then 
$$z\equiv gt+t^{-1} \ , \ x_2\equiv (\epsilon\beta)^{-1}(gt-t^{-1}).$$
Thus 
$$x_1\equiv 2^{-1}(\pm z-x_2)\equiv 2^{-1}(\pm (gt+t^{-1})-(\epsilon\beta)^{-1}(gt-t^{-1})).$$
There are three values of $x$ that map to the same $h=x(x^2-g)$, except 
for in possible special cases.

Here $gt-t^{-1}\not\equiv 0$ and $gt+t^{-1}\not\equiv 0$ as $-1\in QR_p$.
Therefore $z$ and $x_2$ are always in $Z_p^*$. It means that 
$x_{1+}\not\equiv x_{1-}$. It is also not possible that $x^2-g\equiv 0$.

The only special case that can appear here is that $x_2\equiv x_{1\pm}$.
If so, then $x_1=2^{-1}(\pm z-x_2)=x_2$. Thus $3x_2=\pm z$. We get the 
equation
$$3(\epsilon\beta)^{-1}(gt-t^{-1})\equiv \pm (gt+t^{-1})$$
which simplifies to
$$3(\epsilon\beta)^{-1}(gt-t^{-1})\equiv \pm (gt+t^{-1}).$$
This leads to
$$t^2\equiv g^{-1}(1\pm 3(\epsilon\beta)^{-1})^{-1}(1\mp 3(\epsilon\beta)^{-1})$$
$$t^2\equiv g^{-1}4(1\pm 3(\epsilon\beta)^{-1})^{-2}.$$
Clearly, there are no solutions $t$ to this equation. 

The result is that if $c=g$ and $3\in QR_p$ (and $-1\in QR_p$ as assumed in
the lemma), then $m=1$ and $m_1=0$.

The proof of Lemma 5 is complete. $\eop$

Let us give two examples of Lemma 5. For $p=13$ we get for the running
index $t=1,2,3,4,5,6,7,8,9,10,11,12$ the following sequences
$x_2=1,-,6,12,3,10,3,10,1,7,-,12$, 
$x_{1+}=0,-,10,0,7,6,3,10,12,3,-,1$,
$x_{1-}=12,-,10,1,3,10,7,6,0,3,-,0$,
$h\equiv x_2(x_2^2-1)=0,-,2,0,11,2,11,2,0,11,-,0$.

The two values of $t$ that do not give $x_2$ are $2$ and $-2=12$.
They are solutions to $3t+t^{-1}\equiv 0$, i.e., 
$t^2\equiv -3^{-1}\equiv -9\equiv 4$, thus $t\equiv \pm 2$.
There are two values of $t$ giving $x_{1+}\equiv x_{1-}$. They are
solutions to $t^2\equiv 3^{-1}\equiv 9$, i.e., $t=3,10$.  
We get these four special values of $t$ because $3\equiv 2^4\in QR_{13}$. 

The four values of $t$ that give $h\equiv x_2(x_2^2-1)$ are values
$t=\pm 1, \pm 4$. They give $x_2\equiv \pm 1$ and thus $x_2^2-1\equiv 0$.
These special values of $t$ appear because $c=1$. 
For $c=g$ there are no such special values of $t$.

There are four values of $t$ giving $x_2\equiv x_{1\pm}$. They are
solutions to $(t\pm 1)^2\equiv 3^{-1}+1\equiv 10$. 
Thus, $t\pm 1=\pm 6$. That yields $t=5,6,7,8$. We get these special values
of $t$ because $3^{-1}+1\equiv 10=2^{10}\in QR_{13}$. 

The sum of the special values of $t$ is $4+4+4=12$. 
This is subtracted from $p-1=12$ and the result is zero
showing that there are no cases when three different values of $x$ map to
the same $h$ for $p=13$.  

As the second example consider $p=29$. There are four groups of triplet
$x$ values that map to the same $h$: 

$x=5,10,14=-24,-19,-15$ map to $4$, and

$x=19,24,15=-10,-5,-14$ map to $25\equiv -4$. While

$x=12,20,26=-17,-9,-3$ map to $5$, and

$x=3,9,17=-26,-20,-12$ map to $24\equiv -5$.

These $12$ values of $x$ need $24$ values of the running index $t$.
Additionally we have two values of $x$ that give $h\equiv 0$, that is,
$x\equiv \pm 1$. These need four values of $t$. This means that
all $p-1=28$ values of $t$ are already used in the case of multiples and
there cannot be any more special values of $t$. This is indeed true since 
$3\in QNR_{29}$ and $3^{-1}+1\equiv 11\in QNR_{29}$. 
Therefore $m=5$, as Lemma 5 says, and the number of
triplets of values of $x$ mapping to the same $h$ is $(p-5)/6=4$. 
Both $h\equiv \pm 5$ are in $QR_{29}$ and both $h\equiv \pm 4$ are in 
$QR_{29}$. The expectation value is that $(p-5)/24=1$ are in $QR_{29}$.
We got one more in $p=29$: there is random variation around the mean.  

There are exactly $14=(p-1)/2$ singletons:
$$x=2,4,6,7,8,11,18,21,23,24,25,27,13,16$$ 
map to 
$$h=6,2,7,17,11,15,14,18,22,12,27,23,9,20.$$
Of these $14$ singleton $h$ values six are in $QR_{29}$, namely
$6,7,22,23,9,20$. It is about half of $14=(p-1)/2$ numbers in $QR_{13}$.
It is not exactly half, there is random variation around the half.

Let us formulate an assumption that allows us to estimate averages and
variances of the numbers $a_p$.

\th Statistical Assumption.
Let the numbers $h$ be defined as 
$$h\equiv s+s^{-1}\ (\mod p_j) \, \ s=1,\dots, p_j-1$$
and let $n_j$ be the number of $s$ that give $h\in QR_{p_j}$. 
Then $n_j/2$ is binomially distributed with the mean at
$(p_j-1)/2$ and the variance $(p_j-1)/4$. Notice that $s+s^{-1}$
gives each value $h$ exactly twice.

\th Lemma 6.
Let $a_{p_j}=p_j-N_{p_j}$ and $N_{p_j}$ be the number of solutions
to 
$$y^2=x^3+ax+b.\eqno(20)$$ 
The following claims hold assuming that the Statistical Assumption holds.

(i) In Cases 1 and 2 of Lemma 3 the expectation values of $a_{p_j}$ are $2$
and $-2$ respectively when $p_j$ ranges over all primes. For Cases 3 and 4
the expectation value of $a_{p_j}$ is zero. 

(ii) The variance of $a_{p_j}$ is $2p_j$ when $p_j$ ranges over all primes.

(iii) The statistical bound by the standard deviation is: $|a_{p_j}|\le 2p_j^{1\over 2}$.

\proof

By the Statistical Assumption the expectation value of $n_1$ in Lemma 3
for Cases 1 and 2 is $E[n_1]={p-5\over 8}$. Thus
$E[a_p]=p-(p-5+7)=-2$ for Case 1 and $E[a_p]=p-(2p-p+5-7)=2$ for Case 2.
In the Cases 3 and 4 of Lemma 3 $E[n_2]={p-3\over 8}$. Thus
$E[a_p]=p-(p-3+3)=0$ for Case 3 and $E[a_p]=p-(2p-p+3-3)=0$ for Case 4.
 
The expressions of $N_p$ in Lemma 3 are not convenient 
for calculating variances needed in (ii) as Lemma 3 does not count how many
times the running index hits the same $x$. 
We have to use the forms in Lemma 5. 

Lemma 5 lists numbers $m\le 13$ and $m_1\le 6$, 
which are different for different cases of multiple $x$ values. 
These few special cases are not important when estimating the 
variance of $a_p$ when $p$ ranges over all primes. For that reason we assume
that $p$ is so large that $13<<p$ and all values of $t$ in the case
of multiples give a value $h$ that comes from three values of $x$.

Case 1: $b=0$. Lemma 4 shows that $-1\in QNR_{p_j}$ in half of the cases 
of $p_j$. 
By Lemma 1 $a_{p_j}=0$ in these cases. In the remaining half of the cases
$-1\in QR_{pj}$. Exactly half of them are singleton cases of Lemma 5 and 
exactly half are multiplet cases of Lemma 5.

In Lemma 5 is used a running index $t=1,..,p-1$. 
In singleton cases two values of $t$ map to the same $x$. 
Only one $x$ maps to a given $h=x(x^2-c)$, $c=1$ or $c=g$. 
Both $h=x(x^2-c)$ and
$-h$ either belong to $QR_{p_j}$ or both belong to $QNR_{p_j}$. For each
$h\in QR_{p_j}$ there are two values $y$, (i.e., $\pm y$) that give the 
same $y^2\equiv h$. Thus, eight values of $t$ map together as a unit.

In multiplet cases two values of $t$ map to the same $x$. Ignoring the
$m\le 13$ special cases, three values of $x$ map to a given $h$. 
Both $h=x(x^2-c)$ and
$-h$ either belong to $QR_{p_j}$ or both belong to $QNR_{p_j}$. For each
$h\in QR_{p_j}$ there are two values $y$, (i.e., $\pm y$) that give the 
same $y^2\equiv h$. Thus, twenty four values of $t$ map together as a unit.
There is no difference between the cases $c=1$ and $c=g$ except for in the
small numbers $m$ and $m_1$ in Lemma 5. 

By the Statistical Assumption we can treat the situation as trials in a 
binomial distribution. In the binomial distribution there are $n$ trials
and the success probability is $p$ and $q=1-p$. Then the expectation value
is $E[$ number of successes $]=np$ and the variance is $npq$. Here we have
$p_j/k$ independent trials, where for singletons $k=8$ and for multiplets
$k=24$. Each success gives $k$ units and the probability $p=0.5$. Thus,
the average $np$ must give $p_j$ as $a_{p_j}=p_j-N_{p_j}$ should have 
the average at zero, or close to zero, for us to use the binomial distribution
assumption. We have a better match to the binomial distribution if we count 
not $(x,y)$ but $(x,y^2)$ cases. Then the expected number of successes from
$p_j$ trials is $p_j/2=p_j\cdot 0.5=p_j/2$. 
In the binomial distribution the variance
would be $p_j\cdot 0.5\cdot 0.5=p_j/4$ 
as $p=q=0.5$, but in our distribution $k/2$ 
units of success are in a group. The number of independent trials is
$p_j/(k/2)$ and the unit size is multiplied by $k/2$. The effect of this
grouping is that the variance is multiplied by $k/2$, i.e., the variance
$\sigma^2={1\over n}\sum (y_i-\bar y_i)^2$ is multiplied by $a$ if $n$ is
chaged to $n/a$ and $y_i$ changed to $ay_i$.  

The variance of the singleton cases 
is therefore $(p_j/4)\cdot (8/2)=p_j$
and the variance of the multiplet case 
is $(p_j/4)\cdot (24/2)=3p_j$.
Singletons represent one fourth of all cases and multiplets represent 
one fourth of all cases, while $-1\in QNR_{p_j}$ are half of all cases.   
Adding the cases with their probabilities of occurance gives
$$\sigma^2_{y^2}={1\over 2}\cdot 0 + {1\over 4}\cdot p_j + {1\over 4}\cdot 3p_j =p_j.$$
Multiplying the result by $2$ to count $\pm y$ instead of $y^2$ gives
$$\sigma=2p_j.$$

Case $b\not=0$. In this case the addition of $b$ removes the condition
that guaranteed that $a_{p_j}=0$ if $-1\in QNR_{p_j}$. 
For the case $-1\in QR_{p_j}$ this same condition guaranteed that $x$
and $-x$ both give $h$ and $-h$ in $QR_{p_j}$ or both give $h$ and $-h$
in $QNR_{p_j}$. This condition is removed from both cases when $b\not=0$ 
because if $h'_1-b\equiv h$ and
$h'_2-b\equiv -h$, then the relation between $h'_1$ and $h'_2\equiv -h'_1+2b$
does not say anything of $h'_2\in QR_{p_j}$ if $h'_1\in QR_{p_j}$. 

In other ways the case $b\not=0$ does not differ. Having $b\not=0$ does
not even change the special cases of $t$. It only removes the condition
that was discusses above. Exactly half of the cases are singletons and 
exactly half of the cases are multiplets. 
The effect of removing the coupling of $x$ and $-x$ causes
that in the singleton case $k=4$ and in the multiplet case $k=12$.
As a compensation, there is no case of $-1\in QNR_{p_j}$ when there
is zero variation of $a_{p_j}$. The case where $-1\in QNR_{p_j}$ is the same
as the case where $-1\in QR_{p_j}$. These cases give the same variance
and the total variance is
$$\sigma^2_{y^2}={1\over 2}({1\over 2}p_j+{1\over 2}3p_j)+{1\over 2}({1\over 2}p_j+{1\over 2}3p_j)=p_j$$
Counting the solutions $(x,y)$ instead of $(x,y^2)$ gives
$$\sigma=2p_j.$$ 

The statistical bound (iii) is a bound by standard deviation
and it is known as the Hasse bound. The
average variance over all $p_j$ is $2p_j$, 
but if $b=0$, then for the part $-1\in QR_{p_j}$
the variance is $4p_j$. Thus the standard deviation 
for these $p_j$ is $2p_j^{1\over 2}$.

The proof of Lemma 6 is finished. $\eop$

\section{On the zeros of the Taylor series}

\th Lemma 7. 
Let $\phi(s)$ have a Taylor series at $s_0$ of the form
$$\phi(s)=C_0+C_r(s-s_0)^r+C_{r+1}(s-s_0)^{r+1}\cdots$$
Let
$$h(s)={d\over ds}\ln \phi(s).$$
Then 
$$h(s)={rC_r\over C_0}(s-s_0)^{r-1}+O((s-s_0)^r)\hskip 2em {\rm if}\hskip 1em C_0\not=0$$
$$h(s)={r\over s-s_0}+O(1)\hskip 2em {\rm if}\hskip 1em C_0=0.$$

\proof
The claim comes directly from calculating
$$\phi'(s)=rC_r(s-s_0)^{r-1}+(r-1)C_{r+1}(s-s_0)^r\cdots$$
$$h(s)={d\over ds}\ln \phi(s)={\phi'(s)\over \phi(s)}.$$
$\eop$

Why Lemma 7 is written down here is that it may not be obvious that
the only singularity $h(s)$ can have at $s=s_0$ is 
a first order pole. The algebraic rank of $\phi(s)$
is a multiplier (i.e., residue) in the divergent part of $h(s)$.
If the residue of $h(s)$  negative, then the function $\phi(s)$ has
a pole of the order ot the residue.

\th Lemma 8. 
Consider an infinite product
$$\phi(s)=\prod_{j\in A}(1-f_j(s))^{-1}\eqno(21)$$
where $A$ is an infinite subset of $\N$.
Let us assume there is $\alpha>0$ and $C>0$ such that 
$$|f_j(s)|<Cp_j^{-\alpha x}\eqno(22)$$
and
$$|f'_j(s)|<C\ln(p_j)p_j^{-\alpha x}.\eqno(23)$$ 
Here $s=x+iy$ and $f'_j(s)$ denotes the derivative of $f_j(s)$.
The following claims hold:

(i) The function $\phi(s)$ is finite and nonzero 
if $Re\{s\}>{1\over 2\alpha}$.

(ii) Let $s_0$ satisfy ${1\over 3\alpha}<Re\{s_0\}\le {1\over 2\alpha}$.
If the function
$$g_1(s)=\sum_{j\in A}f'_j(s)f_j(s)\eqno(24)$$
is finite at $s_0$, then 
$$h(s)=\sum_{j\in A}f'_j(s).\eqno(25)$$
is finite at $s_0$ and $\phi(s_0)$ is finite and nonzero.

\proof
From the assumption (22) follows that the infinite product (21) 
is absolutely convergent if $Re\{s\}=x>\alpha^{-1}$. 
We expand
$${d\over ds}\ln \phi(s)=\sum_{j\in A}\left(f'_j(s)\sum_{k=0}^\infty f_j(s)^k\right).\eqno(26)$$
By (22) $|f_j(s)|<Cp_j^{\alpha x}$. 
Because
$$\sum_{j=1}^\infty p_j^{-x}$$
is absolutely convergent if $x>1$, the derivative of the series
$$\sum_{j=1}^\infty \ln (p_j) p_j^{-x}$$
is absolutely convergent if $x>1$. By (23) holds
$$|f'_j(s)|<C\ln(p_j)p_j^{-\alpha x}.$$
Thus the series 
$$\sum_{j\in A}f'_j(s)$$
is absolutely convergent if $Re\{s\}>\alpha^{-1}$. 

In the right side series of (26) the terms for each $k$ converge absolutely
if $Re\{s\}>{1\over (k+1)\alpha}$. 
The terms for $k>0$ define an analytic function
$$g(s)=\sum_{j\in A}\left(f'_j(s)\sum_{k=1}^\infty f_j(s)^k\right)$$
This function is finite if $Re\{s\}>{1\over 2\alpha}$. The function
$$g_1(s)=\sum_{j\in A}f'_j(s)f_j(s)$$
in (11) is the first part of this series. The other parts converge if
$Re\{s\}>{1\over 3\alpha}$. 
The function $g_1(s)$ can be analytically continued to the area
${1\over 2\alpha}>Re\{s\}\ge {1\over 3\alpha}$. 
The continuation is 
analytic with the exception of possible isolated singularities, poles.  
The series for $k=0$ converges if $Re\{s\}>{1\over \alpha}$ and defines
an analytic function 
$$h(s)=\sum_{j\in A}f'_j(s).$$
This function can be analytically continued to 
$Re\{s\}>{1\over 3\alpha}$. 
The continuation is 
analytic with the exception of possible isolated singularities, poles.  

Whether $h(s)$ or $g_1(s)$ 
is infinite in a given point $s$ or not, we can formally write
$${d\over ds}\ln \phi(s)=\phi(s)^{-1}{d\over ds}\phi(s)=h(s)+g(s)$$
$${d\over ds}\phi(s)=\phi'(s)=h(s)\phi(s)+g(s)\phi(s)$$
Assume $\phi(s_0)=0$ for $s_0=x_0+iy_0$, 
${1\over \alpha}>x_0\ge {1\over 2\alpha}$.
The function $g(s)$ is finite at $s_0$, thus $g(s_0)\phi(s_0)=0$. If
$h(s)$ is finite at $s_0$, then $h(s_0)\phi(s_0)=0$ and $\phi'(s_0)=0$.
All derivatives of $g(s)$ are finite at $s_0$. 

If $h(s)$ is finite at
$s_0$, then all derivatives of $h(s)$ are finite at $s_0$. If so, 
then by induction all derivatives of $\phi(s)$ are zero at $s_0$:
assuming that it is proved that ${d^j\over ds^j}\phi(s)$ for $j<n$, 
then ${d^n\over ds^n}\phi(s)$ is given by a sum where every term is
a finite value from derivatives of $g(s)$ or $h(s)$ multiplied by
zero, a derivative $j\ge 0$ of $\phi(s)$ at $s_0$. If all derivatives 
of $\phi(s)$ are zero at $s_0$, then $\phi(s)$ is zero everywhere. This
is a contradiction as $\phi(s)$ is described by an infinite product (21)
that converges and is not zero when $Rs\{s\}>\alpha^{-1}$. 
It follows that $\phi(s_0)$ cannot be zero if $h(s_0)$ is finite. 

Next, consider the infinite product of two complex variables $s$ and $z$
$$\psi(s,z)=\prod_{j\in A}(1-f_j(s)+f_j(s+z))^{-1}.\eqno(27)$$
The infinite product is absolutely convergent if $Re\{s\}>\alpha^{-1}$. 
Let us expand 
$${\partial\over \partial z}\ln \psi(s,z)=-\sum_{j\in A}\left(f'_j(s+z)\sum_{k=0}^\infty (f_j(s)-f_j(s+z))^k\right)\eqno(28)$$
The terms $k>0$ define a function of two complex variables
$$u(s,z)=-\sum_{j\in A}\left(f'_j(s+z)\sum_{k=1}^\infty (f_j(s)-f_j(s+z))^k\right)$$
This function is finite if $Re\{s\}>{1\over 2\alpha}$.

The series for $k=0$ converges absolutely 
if $Re\{s\}>{1\over \alpha}$ 
and defines a complex analytic function 
$$h(s+z)=\sum_{j\in A}f'_j(s+z).$$
This function is the same function $h(s)$ as earlier. Now it only has
the variable $s+z$ instead of $s$. We assume it is already continued 
analytically to $Re\{s\}>{1\over 3\alpha}$. It may have isolated poles
when $Re\{s\}>{1\over \alpha}$.

We can write
$${\partial\over \partial z}\ln \psi(s,z)=\psi(s,z)^{-1}{\partial\over \partial s}\psi(s,z)=h(s+z)+u(s,z)$$
$${\partial\over \partial z}\psi(s,z)=h(s+z)\psi(s,z)+u(s,z)\psi(s,z)\eqno(29)$$
The functions $g(s)$ and $u(s,z)$ are analytic (implying finite)
at an open environment of the point $s=s_0$ 
because ${1\over \alpha}\ge Re\{s\}>{1\over 2\alpha}$. 

The function $\psi(s,s_0-s)$ is an analytic function of one complex variable
$s$ and $\psi(s_0,0)=1$. 
The function $\psi(s,s_0-s)$ cannot be zero in every point in an open 
environment of $s_0$ because else it is zero everywhere.
Likewise, $\psi(s,s_0-s)$ cannot be
infinite in every point in an open environment of $s_0$ because then 
$\psi(s,s_0-s)^{-1}$ is zero everywhere. 
It follows that we find $z\not= 0$ such
that $|z|<<1$, $s_0=s+z$ and $\psi(s,s_0-s)\not= 0$. Let this point be 
$(s_2,z_2)$. It follows that $s_2+z_2=s_0$.
Keeping $s$ fixed at $s_2$ and letting $z$ vary, $\psi(s_2,z)$ defines
a complex analytic function of $z$. As $\psi(s_2,z_2)$ is finite, the 
partial derivative of $\psi(s_2,z)$ with respect to $z$ at $z=z_2$ is finite.
Thus, the left side of (29) is finite. The right side contains a finite term
$u(s_2,z_2)\psi(s_2,z_2)$ and $\psi(s_2,z_2)$ is finite and nonzero.
It follows that $h(s_2+z_2)=h(s_0)$ is finite. 
Because $h(s_0)$ is finite $\phi(s_0)$ is nonzero. The first claim (i) of 
the lemma is proven.

In the second claim (ii) of the lemma we select 
a point $s=s_0$ such that ${1\over 2\alpha}\ge Re\{s\}>{1\over 3\alpha}$.
The function $g(s)-g_1(s)$
is finite in this area because the series defining this function
converges absolutely.
The function $u(s,z)$ converges with the exception of the first term. Thus 
$$u_2(s,z)=-\sum_{j\in A}\left(f'_j(s+z)\sum_{k=2}^\infty (f_j(s)-f_j(s+z))^k\right)$$
is finite in this area. In the first term 
$$u_1(s,z)=-\sum_{j\in A}f'_j(s+z)(f_j(s)-f_j(s+z))$$
we insert $z=s_0-s$ and get a complex analytic function of one variable $s$.
This function can be continued analytically to the area
${1\over 2\alpha}\ge Re\{s\}>{1\over 3\alpha}$.
If $g_1(s_0)$ is finite, then also $u_1(s_0,s_0-s_0)=u_1(s_0,0)$ is
finite. The argument that $h(s_0)$ is finite goes in the same way as 
in the case where 
${1\over 2\alpha}< Re\{s\}\le {1\over \alpha}$. 
This proves the second claim.
$\eop$

\section{Application to $L(C,s)$}

\th Lemma 9.
Let $a_{p_j}=p_j-N_{p_j}$ 
and let $N_p$ be the number of solutions to 
$$y^2\equiv x^3-d^2x\hskip 1em (\mod p).$$

(i) The function $g_1(s)$ with the series expression
$$g_1(s)=\sum_{j\in A}\ln(p_j)a_{p_j}^2p_j^{-2s}$$
in the area where the series converges does not depend 
on $d$.

(ii) If $d=k^2$ the function $h(s)$ with the 
series expression 
$$h(s)=\sum_{j\in A}\ln(p_j)a_{p_j}p_j^{-s}$$
in the area where the series converges does not depend 
on $k$. 

\proof
From Lemma 1 follows that if $-1\in QNR_{p_j}$ then $a_{p_j}=0$. 
If $-1\in QR_{p_j}$ Lemma 2 gives two possible values for $a_{p_j}$:
$$a_{p_j}=p_j-N_{p_j}=p_j-8n_{1,j}-7 \ {\rm if \ } d\in QR_{p_j}$$
$$a_{p_j}=p_j-N_{p_j}=p_j+8n_{1,j}+7 \ {\rm if \ } d\in QNR_{p_j}$$
where $n_{1,j}$ is the number of solutions $y^4\in Z^*_{p_j}$ 
yielding $y^4-1\in QR_{p_j}$. The number $n_{1,j}$ does not depend on $d$,
thus the square of either of $s_{p_j}$ does not depend on $d$. It
follows that $g_1(s)$ does not depend on $d$.

In $h(s)$ the number $d$ is a square $k^2$. 
Therefore $d$ is a square in every $Z^*_{p_j}$. 
It follows that for every $p_j$ the case in Lemma 3 is
always Case 1. The value $N_{p_j}=8n_{1,j}+7$ 
for every $p_j$. Thus, the function $h_0(s)$ does not depend
on what square number $d=k^2$ is used: in Lemma 2 the number $a=-d^2$ is 
removed at the beginning by a substitution and $c=1$ in Case 1 of Lemma 2.
$\eop$

The problem statement [1] says that the product (1) converges 
absolutely for $Re\{s\}>{3\over 2}$. This claim assumes the 
statistical Hasse bound $|a_{p_j}|<2p_j^{1\over 2}$. Therefore
this bound must be considered accepted in the context of the problem
statement of [1]. The Hasse bound follows from the Statistical Assumption
and implies that some similar statistical assumption is used in the bound.

\th Lemma 10.
If $\phi(s)$ is Lemma 8 is $L(C,s)$ as defined in (1), then
the term of the infinite product in Lemma 8 is
$$f_j(s)=a_{p_j}p_j^{-s}-p_j^{1-2s}.$$
Then
$$h(s)=\sum_{j\in A}f'_j(s)=h_1(s)+h_2(s)$$
where
$$h_1(s)=-\sum_{j\in A}a_{p_j}\ln(p_j)p_j^{-s}.$$
The function $h_2(s)$ is
$$h_2(s)=2\sum_{j\in A}\ln(p_j)p_j^{1-2s}.$$
It diverges at $s=1$ and has a simple pole
$$h_2(s)={1\over s-1}+{\rm finite} \ {\rm terms}$$
at $s=1$.
The function $g_1(s)$ we define as
$$g_1(s)=\sum_{j\in A}f'_j(s)f_j(s)=g_{1,1}(s)+g_{1,2}(s)$$
where
$$g_{1,1}(s)=\sum_{j\in A}\ln(p_j)a_{p_j}^2p_j^{-2s}.$$
It is the the part of $g(s)$ that diverges if $s=1$ and $g_{1,2}(s)$
and the other parts of $g(s)$ converge at $s=1$.

\proof
A simple calculation shows that the forms of $h(s)$ and $g_1(s)$ are
as in the claim of the lemma. 
   
The pole of the function $h_2(s)$ at $s=1$ is derived from the pole
of the Riemann zeta function $\zeta(s)$ at $s=1$. Zeta has a simple pole
of residue $1$ at $s=1$. Thus close to $s=1$ zeta is
$$\zeta(s)=\prod_{p_j}(1-p_j^{-s})^{-1}={1\over s-1}+{\rm finite} \ {\rm terms}.$$
Derivating gives $h(s)$ for the zeta function:
$${\zeta'(s)\over \zeta(s)}=-{1\over (s-1)^2}\cdot {s-1\over 1}+{\rm finite} \ {\rm terms}=h(s)+g(s).$$
Calculating from the infinite product we get
$${d\over ds}\ln(\zeta(s))=h(s)+g(s)$$
where
$$h(s)=-\sum_j\ln(p_j)p_j^{-s}$$
$$g(s)= -\sum_j\ln(p_j)p_j^{-2s}+\cdots$$
and $g(s)$ converges when $Re\{s\}>{1\over 2}$. Thus, the divergent 
part is
$$\sum_j\ln(p_j)p_j^{-s}={1\over s-1}+{\rm finite} \ {\rm terms}$$
close to $s=1$. Changing $-s=1-2z$ gives 
$$2\sum_j\ln(p_j)p_j^{1-2z}=2{1\over 2z-2}+{\rm finite} \ {\rm terms}$$
that is
$$2\sum_j\ln(p_j)p_j^{1-2s}={1\over s-1}+{\rm finite} \ {\rm terms}.$$

In Lemma 6 we derived the Hasse bound $|a_{p_j}|\le 2p_j^{1\over 2}$. 
It is a statistical bound where a random variable is estimated by 
the standard deviation. This statistical bound can be used for 
investigation of the convergence of the $L(C,s)$ function because 
any looser bound 
$|a_{p_j}|\le 2p_j^{{1\over 2}+\alpha}$, $\alpha>0$, on a general $p_j$
can be improved by adding the probability that $|a_{p_j}|$ is higher than 
(a chosen multiple) of the standard deviation of $|a_{p_j}|$.

Lemma 8 shows that the bound on $<a_{p_j}|$ cannot be a
stricter bound of the type $|a_{p_j}|\le p_j^{{1\over 2}-\alpha}$
where $\alpha>0$ for a general $p_j$ 
because if there is such a bound, then $g(s)$
converges at $s=1$ and Lemma 8 shows that $L(C,s)$ is always finite
and nonzero.
Yet we know that $\ln(L(C,s))$ is infinite at $s=1$ for some elliptic
curves $C$. 

From the Hasse bound it follows that $g_1(s)$ converges for $Re\{s\}>1$.
By Lemma 8 $g_1(s)$ must diverge at least for some $C$ at $s=1$
because there are elliptic curves $C$ such that $L(C,1)$ is not finite and
nonzero. The Hasse bound guarantees that
$g(s)-g_1(s)$ converge absolutely if $Re\{s\}>{1\over 2}$ and 
that $g_{1,2}(s)$ also converges absolutely if $Re\{s\}>{1\over 2}$.
The part that cannot be shown to converge is $g_{1,1}(s)$.
It must diverge at $s=1$ because $g_1(s)$ is the same function
for all elliptic curves $C$ of the type $y^2=x(x^2-d^2)$ (see Lemma 9) and 
for some values of $d$ the L-function has a zero value at $s=1$.
If $g_{1,1}(s)$ would converge at $s=1$, then $L(C,1)$ would
be finite and nonzero for every $d$.

All claims of the lemma are proven.$\eop$

With these results we can discuss the Birch and Swinnerton-Dyer conjecture. 
The original reason why Birch and Swinnerton-Dyer formulated their conjecture
was that the function
$$\log \prod_{j=1}^n{N_{p_j}\over p_j}$$
was growing approximatively linearly for some values of $d$ for
elliptic curves of the type $y^2=x(x^2-d^2)$, while for some other 
values of $d$ the function did not tend to infinity when $n$ grows.

The function they studied can be written as
$$\log \prod_{j=1}^n{N_{p_j}\over p_j}=-\log \prod_{j=1}^n(1-{a_{p_j}\over p_j})^{-1}$$
$$=\sum_{j=1}^n\log(1-a_{p_j}p_j^{-1})=-\sum_{j=1}^n a_{p_j}p_j^{-1}+\sum_{j=1}^n a_{p_j}^2p_j^{-2}+\cdots$$
As $a_{p_j}$ is on the range of $p_j^{1\over 2}$ by the Hasse bound, the
higher terms can be ignored. The first two terms cannot be ignored.
The derivative of the first one is $h(1)$ and $h(s)$ can diverge. 
The derivative of the second one is $g_{1,1}(1)$
and Lemma 8 shows that it diverges at least if $b=0$ since 
$L(C,1)$ is not nonzero for all elliptic curves with $a=-d^2$, $b=0$. 

The function Birch and Swinnerton-Dyer studied can be completed to a 
complex function
$$L_2(C,s)=\prod_{j\in A}(1-a_{p_j}p_j^{-s})^{-1}.$$
The difference with this function and $L(C,s)$ in (1) is that the
term $p_j^{1-2s}$ is missing. 
For $L_2(C,s)$ the functions in Lemma 10 are
$$h(s)=\sum_{j\in A}f'_j(s)=h_1(s)$$
where
$$h_1(s)=-\sum_{j\in A}a_{p_j}\ln(p_j)p_j^{-s}.$$
The divergent part of the function $g_1(s)$ is 
$$g_{1,1}(s)=\sum_{j\in A}\ln(p_j)a_{p_j}^2p_j^{-2s}.$$
Thus, $g_{1,1}(s)$ is the same for $L_2(C,s)$ and for $L(C,s)$ but
the function $h(s)$ lacks the second part $h_s(s)$ in Lemma 10. This
this missing function has a first order pole and residue one at $s=1$.

\th Theorem 1. The Birch and Swinnerton-Dyer conjecture, as formulated
in [1], fails for rank zero.

\proof
Birch and Swinnerton-Dyer studied elliptic curves of the form 
$y^2=x(x^2-d^2)$ with several values of $d$. As these curves have rank 
zero if $d=k^2$, it must be so that the function they studied, $L_2(C,s)$,
has a finite nonzero value at $s=1$ for $d=k^2$. Consequently
$$h(s)+g_{1,1}(s)$$ 
must be finite at $s=1$. As $L_2(C,s)$ does not have the $h_2(s)$ part, 
$h(s)=h_1(s)$. 

Consequently, the function $L(C,s)$ cannot have a finite nonzero value
at $s=1$ for elliptic curves with $a=-k^4$ and $b=0$. The function $h_2(s)$
has a pole at $s=1$ and is of the form
$$h_2(s)={1\over s-1}+{\rm finite} \ {\rm terms}.$$
by Lemma 9. 
This result does not need any statistical assumptions.
It is a direct result of the pole of the Riemann zeta at $s=1$.
Therefore
$$h_1(s)+h_2(s)+g_{1,1}(s)={1\over s-1}+{\rm finite} \ {\rm terms}.$$ 
Because of this pole the function $L(C,s)$ has a zero at $s=1$. Thus, 
for an elliptic curve of rank zero the L-function $L(C,s)$ has a zero
at $s=1$ and has the algebraic rank one.$\eop$

The Theorem on page 4 in [1] claims that it is proven that if the elliptic
curve $C$ has rank zero then the L-function $L(C,s)$ has the algebraic
rank zero. This theorem is false. It is in direct contradiction with the
initial experiments of Birch and Swinnerton-Dyer.

We can prove that the initial experiments of Birch and Swinnerton-Dyer were
correct. The function $L_2(C,s)$ does indeed have a nonzero value at $s=1$
if $C$ is an elliptic curve with $a=-k^4$ and $b=0$. 

The function $L_2(C,s)$ necessarily has a pole since Birch and
Swinnerton-Dyer found that if $d=5$ the logarithm of $L_2(C,s)$ grows
to infinity. A logarithm grows to infinity in two cases: either the
function $L_2(C,s)$ has a zero at $s=1$ or it has a pole at $s=1$.
It must be the second case since for $d=k^2$ the function $h_1(s)$ 
grows to positive infinity at $s=1$, as is shown in Lemma 11. When $d=k^2$
every $a_{p_j}$ is taken form Case 1 of Lemma 3 and the expected value
of $a_{p_j}$ is $-2$. In the case $d=5$ about half of the terms
$a_{p_j}$ come from Case 2 of Lemma 3 and they have the expectation value
$2$. Therefore the function $h_1(s)$ for $d=5$ must grow slower than
for $d=k^2$ and so much slower than for $d=5$ the function $\ln(L_2(C,s))$
grows to infinity. It follows that for $d=5$ the function $L_2(C,s)$ has
a pole at $s=1$.

It is proven that $L(C,s)$ is
analytic in the whole plane, see [1]. Thus, it does not have poles.
The pole of $L_2(C,s)$ is removed by the addition of the term
$p_j^{1-2s}$ to the Euler product. This addition gives the function
$h_2(s)$. The function $h_2(s)$ has a simple pole with residue 1 
at $s=1$. As this function removes the pole of $g_{1,1}(s)$ and 
$g_{1,1}(s)$ does not depend on $d$ by Lemma 9, it follows that
$$g_{1,1}(s)={-1\over s-1}+{\rm finite} \ {\rm terms}.$$
Since $h_1(s)+g_{1,1}(s)$ must be a finite function if $d=k^2$ by
the numerical experiments by Birch and Swinnerton-Dyes, 
then $h_1(s)$ must have the form
$$h_1(s)={1\over s-1}+{\rm finite} \ {\rm terms}.$$

The Birch and Swinnerton-Dyer conjecture also does not hold for the function
$L_2(C,s)$.
The function $\log L_2(C,s)$ does grow infinitely for an elliptic curve
with $a=-d^2$, $b=0$, $d=5$, and stays finite if $d=k^2$, but the function
$L_2(C,s)$ does not have a zero when the elliptic curve with $d=5$ has
rank one. The function $L_2(C,s)$ has a pole at $s=1$. 

We can make the following two observations. 

(1) We showed that if $a=-d^2$ and $b=0$ the function
$g_{1,1}(s)$ has at $s=1$ a simple pole of the form
$$g_{1,1}(s)={-1\over s-1}+{\rm finite} \ {\rm terms}\eqno(30)$$
If we replace $a_{p_j}^2$ in $g_{1,1}(s)$ 
$$g_{1,1}=-\sum_{j\in A}{a^2_{p_j}\over p_j}\ln(p_j)p_j^{1-2s}$$
by its average value $2p_j$ from Lemma 6 we get the approximation
$$g_{1,1}=-\sum_{j\in A}2\ln(p_j)p_j^{1-2s}={-1\over s-1}+{\rm finite} \ {\rm terms}$$

(2) We  also showed that for $d=k^2$ the function $h(s)$ with the 
series expression 
$$h(s)=\sum_{j\in A}\ln(p_j)a_{p_j}p_j^{-s}$$
in the area where the series converges is of the form
$$h_1(s)={1\over s-1}+{\rm finite} \ {\rm terms}\eqno(31)$$
when $s$ is close to $1$.
If $-1\in QNR_{p_j}$ the value of $a_{p_j}$ is zero by Lemma 1.  
If for $-1\in QR_p$ we replace $a_{p_j}$ by its average value $-2$ from 
Lemma 6 we get the approximation
$$h(s)={1\over 2}\sum_{j\in A}2\ln(p_j)p_j^{-s}$$
which is exactly of the correct form (31).

Both observations support the idea that the Statistical Assumption
can be used in estimating the values of $a_{p_j}$. 
If the Statistical Assumption holds, the function $g_{1,1}(s)$ 
has a simple pole of the same form (30) for any elliptic
curve. The function $h_1(s)$ also cannot be different: either it is 
finite or it has the form (31). Therefore it is
not possible that the L-function can give algebraic ranks of
six or so on, though elliptic curves can have such ranks.

\section{ Discussion}

I wrote the original version [2] of this paper in 2008. In that paper I had
a simple proof that the Birch and Swinnerton-Dyer conjecture fails for
elliptic curves of the type $y^2=x^3-d^2x$, that is, the type used 
by Birch and Swinnerton-Dyer in their numerical calculations leading to
the conjecture.

I took three values of $d$ for this special type of elliptic curve. 
One value was $d=1$, which gives the rank zero and the other two values 
of $d$ were $d=5$, which gives rank one, and $d=19$, which gives rank zero. 

As my proof showed that the Theorem on page 4 in the problem statement [1] 
is wrong, and as that theorem was based on results in peer-reviewed and 
accepted articles, I did not want to use a literature result that $19$ is
a noncongruent number: which literature result would then be the more 
reliable? Thus, I proved that $19$ is a noncongruent number. The lemmas
in [2] proving this fact are now moved to [3], which I see no reason to 
publish, but I do see reason to check published results. A similar case 
happened to me already before with another Millennium Prize problem.
A literature result was in contradiction with my elementary calculations
in a solution of the Navier-Stokes problem. Referees of two journals used this
conflict for rejecting my paper, but it turned out that the literature 
result, also implicitly used by the problem statement, had an error. My 
elementary calculations were quite correct. 

The simple proof in the original version of [2] was to take the function
$$\ln (L(C_1,s)L(C_2,s)^{-1})=-H_{d_1}(s)+H_{d_2}(s)$$
where
$$H_d(s)=-\sum_{j\in A}a_{p_j,d}p_j^{-s},$$ 
$C_1$ is an elliptic curve with $a=-d_1=-1$, $b=0$ and $C_2$ is 
another elliptic curve where $a=-d^2$, $b=0$. Alternatively 
one can select $d_2=5$ or $d_2=19$. 

The value $d_1=1$ is a square number and picks up only Case 1 from Lemma 3. 
In Case 1 the mean of $a_{p_j}$ is $-2$. This does not require a statistical
assumption of the distribution, as the Hasse bound does. The mean is $-2$
with the simple assumption that the mean of $n_1$ in Lemma 3 is zero.
The value $-2$ comes from the special values of $x$ in Lemma 3 in Case 1, 
i.e., that $x^2\equiv 1$ has two solutions. 
It is always the case as $1\in QR_p$. Thus,
no special assumptions are needed for concluding that $H_1(s)$ does
diverge at $s=1$ and $H_1(s)$ tends to the positive infinite. Whatever
other value $d$ is chosen for the other curve $C_2$, it always has
both Case 1 and Case 2 in Lemma 3, basically it has equally many of
both cases. If $d_2\in QR_p$, then $a_{p_j,d}=a_{p_j,1}$, while
$d_2\in QNR_p$, then $a_{p_j,d}=-a_{p_j,1}$. As there are (about) equally
many cases $d\in QR_p$ as cases with $d\in QNR_p$ for $d=5$ or $d=19$
we get the function 
$$G(s)=\sum_{j\in A}c_jp_j^{-s}H_1(s)-H_d(s)$$ 
where $3/4$ of $c_j$ are zero and the remaining $1/4$ are $2a_{p_j,1}$.
Half of $a_{p_j,1}$ are zero for any elliptic curve with $b=0$
by Lemma 1, thus $G(s)$ must have the same asymptotic behavior as 
$H_1(s)$ and therefore $H_2(s)$ must be finite at $s=1$. This must be
the same for both $d=5$ and $d=19$, thus the Birch and Swinnerton-Dyer
conjecture fails for ranks zero and one and Theorem on page 4 in [1] fails.

I sent my proof to one expert. He discarded it because the result was in
contradiction with what he believed was proven. He did not show errors
in my calculations and indeed, my calculations were correct. It was just
a question of whether to trust authorities or check elementary calculations.
Another expert complained 
that I prove that $19$ is noncongruent though it is proven in the literature.
I did indeed reprove a known result, but it was not possible initially to
know what literature result is false and what is correct. I removed the
claim that the paper contains a counterexample to the Birch ad Swinnerton-Dyer
conjecture and put the paper to arXiv in 2008.  

Now it is 2020 and I have checked my calculations. 
They are totally elementary, well checked by examples, and they are correct. 
The original 2008 version had only the lemmas 1-4
and the lemmas in [3] proving that $19$ is a noncongruent number. 
In order to address the 2008 conterarguments to my proof I have added a new
part that completely explains what happens to $L(C,s)$ at $s=1$. The 
central part of this explanation is Lemma 8. The method of the proof of
Lemma 8 is from [4], but the results from [4] are not needed in the
present paper. Using Lemma 8 for $L(C,s)$ requires the Hasse bound
for convergence of the L-function. Though
[1] assumes that this bound is proven in stating that $L(C,s)$ 
converges absolutely if $Rs\{s\}>3/2$, I think that proofs of well-known
problems should always be elementary and self-contained, not relying on
any published results and especially not on results proven by advanced
methods (i.e., anything too difficult for a second year student with a
simple pocket calculator and a faint memory of algebra I). 
Therefore I added Lemmas 5 and 6 for proving this bound. 
The new Lemma 7, Lemmas 9 and 10 are needed for applying Lemma 8
to $L(C,s)$ and finally this discussion is necessary for replying to the
agrument that the conjecture is proven for ranks zero and one.

\end{document}